\documentclass[final,onefignum,onetabnum]{siamart171218}

\usepackage{lipsum}
\usepackage{amsfonts}
\usepackage{graphicx}
\usepackage{epstopdf}
\usepackage{algorithmic}
\ifpdf
  \DeclareGraphicsExtensions{.eps,.pdf,.png,.jpg}
\else
  \DeclareGraphicsExtensions{.eps}
\fi

\newsiamremark{remark}{Remark}
\newsiamremark{hypothesis}{Hypothesis}
\crefname{hypothesis}{Hypothesis}{Hypotheses}
\newsiamthm{claim}{Claim}

\headers{Maximum Entropy Snapshot Sampling for Reduced Basis Generation}{F. Kasolis and M. Clemens}

\title{Maximum Entropy Snapshot Sampling\\ for Reduced Basis Generation\thanks{Submitted to arXiv.org on \today.
\funding{This work was partly funded by the German Research Foundation (DFG) under grant no. CL143/18-1.}}}

\author{Fotios Kasolis\thanks{Chair of Electromagnetic Theory, University of Wuppertal (\email{kasolis@uni-wuppertal.de}).}
   \and Markus Clemens\footnotemark[2]
}

\usepackage{amsopn}

\ifpdf
\hypersetup{
  pdftitle={Maximum Entropy Snapshot Sampling for Reduced Basis Generation},
  pdfauthor={F. Kasolis and M. Clemens}
}
\fi

\begin{document}

\maketitle

\begin{abstract}
Snapshot back-ended reduced basis methods for dynamical systems commonly rely on the singular value decomposition of a matrix whose columns are high-fidelity solution vectors. An alternative basis generation framework is developed here. The advocated maximum entropy snapshot sampling (MESS) identifies the snapshots that encode essential information regarding the system's evolution, by exploiting quantities that are suitable for quantifying a notion of dynamical stability. The maximum entropy snapshot sampling enables a direct reduction of the number of snapshots. A reduced basis is then obtained with any orthonormalization process on the resulting reduced sample of snapshots. The maximum entropy sampling strategy is supported by rigid mathematical foundations, it is computationally efficient, and it is inherently automated and easy to implement.
\end{abstract}



\section{Introduction}

Problem specific reduced bases are often a key ingredient of reduced basis model reduction strategies. The manufactured low-dimensional bases are typically employed within a projection framework, and they enable the reduction of large-scale problems \cite{Quarteroni2016}. In this context, the proper orthogonal decomposition holds a prominent position. In the proper orthogonal decomposition framework, information regarding the temporal or parametric dependence of the original high-fidelity problem is commonly extracted from homogeneously sampled high-fidelity snapshots with the singular value decomposition \cite{Golub2013}, while ad hoc criteria are used for deciding the sample size. In general, the potential of a reduced basis to represent arbitrary states of the original system, depends strongly on the sampled high-fidelity snapshots and on the method that is used for extracting the features of these snapshots.

Typically, a reduced basis consists of a few left singular vectors of the snapshot matrix that are computed with a variant of the singular value decomposition. The relative magnitudes of the singular values are used as a measure of reducibility and to determine the dimension of the reduced basis. The singular value decomposition is an established data analysis factorization, while it constitutes a fundamental ingredient of miscellaneous methods, such as low-rank matrix approximations, the principal component analysis, and the proper orthogonal decomposition, among others. Being inherently linear and optimal in the least-squares sense, the singular value decomposition is robust for data that are obtained from linear problems, while when it is used for generating reduced bases for nonlinear problems, it removes high-frequency components that are typically present and relevant to the evolution of such systems. In the context of nonlinear dynamical systems, variants of the singular value decomposition have been used for estimating the dimension of strange attractors and the generalized alignment index \cite{Hediger1990, Kantz2004}. Further, kernel methods, such as the kernel principal component analysis \cite{Scholkopf1998,Wirtz2014} and the kernel entropy component analysis \cite{Jenssen2010}, have been proposed for model reduction of nonlinear problems and feature extraction.

Here, based on the asymptotic properties of measure-preserving transformations, we develop a discrete reduced basis framework for nonlinear time-dependent problems. The proposed method is based on a variant of the so-called Grassberger-Procaccia correlation sum \cite{Grassberger1983}, on recurrence quantification analysis \cite{Marwan2007}, and on an estimate of the invariant Kolmogorov-Sinai entropy \cite{Mane1987, Cornfeld1989, Takens1998, Downarowicz2011}. This reduction results from constraining an entropy estimate to be a strictly increasing function of the time index. Then, a reduced basis can be generated with any orthonormalization process, such as the orthogonal-triangular decomposition, of the already reduced snapshot matrix.

The main contribution of this work constitutes of a global sampling interruption criterion and a snapshot sampling method that enables a direct reduction of the snapshot matrix. Since the proposed sampling method is inherently recursive, it can be used as part of the high-fidelity solver, and hence, unnecessary data storage is avoided. The effect of calling a basis generation algorithm after selecting a subset of the high-fidelity snapshots has a significant impact on the required computational time, since it is often the case that the complexity of basis generation algorithms is quadratic in the number of snapshots. This feature makes basis generation feasible for problems that are larger than those that can be handled by the singular value decomposition.

The rest of the paper is organized as follows. In \cref{sec:prel}, we summarize the necessary concepts from the theory of dynamical systems and we introduce projection-based reduced basis methods. In \cref{sec:mess}, we motivate an estimate of the Kolmogorov-Sinai entropy and present the main theorems that regard the maximum entropy snapshot sampling. Compression and reduced basis model reduction experiments follow the theoretical presentation, while we complete our study with a brief summary.

\section{Preliminaries}\label{sec:prel}
Before introducing the maximum entropy snapshot sampling, we present established results from the theory of dynamical systems \cite{Cornfeld1989, Downarowicz2011, Katok1995, Mane1987} and from linear algebra \cite{Golub2013,Quarteroni2016}. In particular, our presentation is restricted to recurrences, ergodicity, entropy, and to the singular value decomposition, in the context of the proper orthogonal decomposition. Our goal is to establish terminology and familiarize the reader with notation that is relevant to the maximum entropy snapshot sampling.

\subsection{Dynamical systems}\label{sec:dynsys}
Let $\Omega\neq\emptyset$ be a set, consider a $\sigma$-algebra $\mathfrak{A}$ of subsets of $\Omega$, and a measure $\mu:\mathfrak{A}\rightarrow[0,1]$ such that $\mu(\Omega)=1$. The tuple $(\Omega,\mathfrak{A},\mu)$ is called a probability space. A measurable measure-preserving transformation $f:\Omega\rightarrow\Omega$ is said to be an endomorphism. By attaching an endomorphism $f:\Omega\rightarrow\Omega$ to a probability space $(\Omega,\mathfrak{A},\mu)$, we define a stationary dynamical system $(\Omega,\mathfrak{A},\mu,f)$. If we assume that the semigroup $(\mathbb{N}_0,+)$ of nonnegative integers acts on $\Omega$ by the iterates of $f$, then $f$ is called a one-sided shift, while the dynamical system itself is said to be time-discrete. In what follows, we only consider stationary dynamical systems whose states $x_n\in\Omega$ are specified by $x_{n}=f^{n}(x)$ for all $n\in\mathbb{N}_0$, where $f^0=\mathrm{id}_\Omega$ is the identity map, each function composition $f^{n+1}=f\circ f^{n}$ specifies an iterate of $f$, and the initial state $x\in\Omega$ is known. A point $x\in A$, with $A\in\mathfrak{A}$, is said to be recurrent, if $f^{n}(x)\in A$ for some integer $n\geq 1$. According to Poincar{\'e}'s recurrence theorem \cite{Cornfeld1989}, \cref{thm:poincare}, almost every $x\in A$ returns to $A$ infinitely often. Apparently, the fixed and the periodic points of an endomorphism, when they exist, are recurrent.

\begin{theorem}[Poincar{\'e}]\label{thm:poincare}
Let $(\Omega,\mathfrak{A},\mu,f)$ be a stationary dynamical system. For any $A\in\mathfrak{A}$ whose measure is nonvanishing, there exists an infinite sequence $(n_k)$ of integers, with $n_k\rightarrow\infty$ as $k\rightarrow\infty$, such that $f^{n_k}(x)\in A$ for almost all $x\in A$.
\end{theorem}
Although \cref{thm:poincare} guarantees the existence of recurrent points, it does not establish quantitative information regarding the number of states in a measurable set of $\Omega$. To this end, Birkhoff's ergodic theorem \cite{Mane1987,Katok1995}, \cref{thm:birkhoff}, ensures the existence of the relative asymptotic number of iterates that are in a measurable set of $\Omega$.

\begin{theorem}[Birkhoff]\label{thm:birkhoff}
Let $(\Omega,\mathfrak{A},\mu,f)$ be a stationary dynamical system. If $\varphi\in L^{1}(\Omega)$, then the limit 
\begin{equation}\label{eq:birkhoff}
\tilde{\varphi}(x)=\lim_{n\rightarrow\infty}\frac{1}{n}\sum_{k=0}^{n-1} \varphi(f^{k}(x))
\end{equation}
exists for almost every point $x\in\Omega$, while $\tilde{\varphi}\in L^1(\Omega)$.
\end{theorem}
A particularly valuable instance of \cref{thm:birkhoff} arises for $\varphi=\chi_{A}\in L^1(\Omega)$, where $\chi_{A}$ is the indicator function of a set $A\in\mathfrak{A}$; that is,
\begin{equation}
\chi_A(x) = \begin{cases} 
      1, & \text{if }x\in A,\\
      0, & \text{if }x\not\in A.
   \end{cases}
\end{equation}
By substituting $\varphi$ for $\chi_{A}$ in \eqref{eq:birkhoff}, we obtain
\begin{equation}\label{eq:sojourn}
\tilde{\chi}_A(x)= \lim_{n\rightarrow\infty}\frac{1}{n}\#\left\{k\in\mathbb{Z}_n~\middle|~f^k(x)\in A\right\},
\end{equation}
where $\#$ is the set cardinality operator and $\mathbb{Z}_n=\{0,1,\ldots,n-1\}$. The number $\tilde{\chi}_A(x)$ is sometimes referred to as the average sojourn time spent by $x$ in the measurable set $A$ \cite{Mane1987}. The limiting value $\tilde{\chi}_{A}(x)$ is not constant for a stationary dynamical system, unless the system's endomorphism is ergodic, according to \cref{def:erg}.
\begin{definition}[Ergodicity]\label{def:erg}
Let $(\Omega,\mathfrak{A},\mu,f)$ be a stationary dynamical system. The endomorphism $f$ is said to be ergodic, if for every set $A\in\mathfrak{A}$ such that $f^{-1}(A)=A$ either $\mu(A)=0$, or $\mu(A)=1$.
\end{definition}
\begin{remark}
In \cref{def:erg}, the notation $f^{-1}(A)$ does not imply that $f$ is invertible, but it denotes the preimage $\{ x\in \Omega \mid f(x)\in A\}$ of $f$.
\end{remark}
Under the assumption that $f$ is ergodic, the measure of a set $A\in\mathfrak{A}$ is equal to the relative asymptotic number of the iterates of $f$ that are in $A$, as follows from the strong ergodic theorem; that is,
\begin{equation}\label{eq:measure}
\mu(A)=\tilde{\chi}_A(x)
\end{equation}
for almost every $x\in\Omega$.

\begin{remark}
In what follows, the term ``system'' refers to a stationary dynamical system whose endomorphism is ergodic.
\end{remark}

Let $(\Omega,\mathfrak{A},\mu,f)$ be a system and consider a finite or countable measurable partition $\mathcal{P}$ of $\Omega$. The Shannon entropy of $\mathcal{P}$ is defined by
\begin{equation}
H^{1}_{\mu}(\mathcal{P}) = -\sum_{\substack{P\in\mathcal{P} \\ \mu(P)>0}}\mu(P)\log{\mu(P)}.
\end{equation}
If $H^{1}_{\mu}(\mathcal{P})<\infty$, then the existence of the limit 
\begin{equation}
h^{1}_{\mu}(f,\mathcal{P}) = \lim_{n\rightarrow\infty}\frac{1}{n}H^{1}_{\mu}\left(\bigvee_{k=0}^{n-1}f^{-k}(\mathcal{P})\right),
\end{equation}
where $\bigvee_{k=0}^{n-1}f^{-k}(\mathcal{P})$ is the joint partition, is guaranteed from the subadditivity of $H^{1}_{\mu}$. Provided the existence of the limit $h^{1}_{\mu}(f,\mathcal{P})$, \cref{def:ksen} enables a metric invariant characterization of a system, which is connected to the system's ability to generate exponentially unstable sequences of states, under perturbations of initial conditions.
\begin{definition}[Kolmogorov-Sinai entropy]\label{def:ksen}
The Kolmogorov-Sinai entropy of a system $(\Omega,\mathfrak{A},\mu,f)$ is
\begin{equation}
\mathcal{H}_{\mu}^{1}(f) = \sup\{h^{1}_{\mu}(f,\mathcal{P})\mid H^{1}_{\mu}(\mathcal{P})<\infty\},
\end{equation}
where the supremum is taken over all measurable partitions $\mathcal{P}$ of $\Omega$.
\end{definition}
For a system $(\Omega,\mathfrak{A},\mu,f)$ whose Kolmogorov-Sinai entropy is (strictly) positive it has been proved \cite{Takens1998} that the Shannon entropy $H^{1}_{\mu}$ can be replaced by the R{\'e}nyi entropy
\begin{equation}
H^{\alpha}_{\mu}(\mathcal{P}) = -\frac{1}{\alpha-1}\log\sum_{\substack{P\in\mathcal{P} \\ \mu(P)>0}}\mu^\alpha(P)
\end{equation}
of any order $\alpha>1$ and the resulting metric invariant R{\'e}nyi entropy $\mathcal{H}^{\alpha}_{\mu}(f)$ is equal to the Kolmogorov-Sinai entropy.

Here, an estimate of the R{\'e}nyi entropy of order $\alpha=2$, which is easy to compute in a numerical setting, turns out to be particularly valuable for efficiently sampling suitable snapshots, in the context of reduced basis generation.

\subsection{The proper orthogonal decomposition}
Let $m$ and $n$ be positive integers such that $m\gg n > 1$, and consider the matrix 
\begin{equation}
\mathbf{X}=(\mathbf{x}_1\mid \mathbf{x}_2\mid\ldots\mid \mathbf{x}_{n})\in\mathbb{R}^{m\times n}
\end{equation}
whose columns $\mathbf{x}_{1}$, $\mathbf{x}_{2}$, $\ldots$, $\mathbf{x}_{n}\in\mathbb{R}^m$ are the iterates of an ergodic endomorphism. In the context of reduced basis methods, the columns of $\mathbf{X}$ are called snapshots, while $\mathbf{X}$ itself is said to be a snapshot matrix. 
\begin{remark}
Although the conditions $m \gg n > 1$ are not necessary, they reflect the fact that we are interested in the reduction of high-dimensional dynamical systems for which it is feasible to obtain more than a single snapshot.
\end{remark}
An orthonormal basis for the range of the snapshot matrix consists of its left singular vectors, and hence, a reduced basis can be found with the singular value decomposition, which is known to exist for all matrices. In the proper orthogonal decomposition framework, $\mathbf{X}$ is factorized according to the singular value decomposition as 
\begin{equation}
\mathbf{X} = \mathbf{U}\mathbf{\Sigma}\mathbf{V}^\top,
\end{equation}
where $\mathbf{U}\in\mathbb{R}^{m\times m}$ and $\mathbf{V}\in\mathbb{R}^{n\times n}$ are orthogonal matrices whose columns are the left and right singular vectors of $\mathbf{X}$, respectively, and $\mathbf{\Sigma}\in\mathbb{R}^{m\times n}$ is the diagonal matrix whose diagonal entries $\sigma_{1}\geq\sigma_{2}\geq\cdots\geq\sigma_{n}\geq 0$ are the singular values of $\mathbf{X}$. Provided $\mathbf{\Sigma}$ and $\epsilon\in(0,1)$, we determine the smallest integer $\ell<\operatorname{rank}(\mathbf{X})$ such that 
\begin{equation}\label{eq:svdinfo}
\Vert\mathbf{\Sigma}(1:\ell,1:\ell)\Vert^{2}_{\mathrm{F}} \geq (1-\epsilon^2) \Vert\mathbf{\Sigma}\Vert^{2}_{\mathrm{F}},
\end{equation}
where the colon slicing notation is described by G. H. Golub and C. F. V. Loan \cite{Golub2013} and $\Vert\ast\Vert_{\mathrm{F}}$ is the Frobenius norm. Then, the first $\ell$ columns of $\mathbf{U}$ constitute a reduced basis for the range of $\mathbf{X}$, while, according to the Eckart-Young-Mirsky theorem, the $\ell$-rank approximation
\begin{equation}
\mathbf{X}_{\ell} = \mathbf{U}(:,1:\ell) \mathbf{\Sigma}(1:\ell,1:\ell)\mathbf{V}(:,1:\ell)^\top
\end{equation}
of $\mathbf{X}$ satisfies the bound $\Vert \mathbf{X} - \mathbf{X}_{\ell} \Vert_{\mathrm{F}} \leq \epsilon\Vert \mathbf{X} \Vert_{\mathrm{F}}$. The matrix 
\begin{equation}\label{eq:projector}
\mathbf{\Pi}_{\ell} = \mathbf{U}(:,1:\ell)\mathbf{U}(:,1:\ell)^\top\in\mathbb{R}^{m\times m}
\end{equation}
is an orthogonal projector on an $\ell$-dimensional subspace of the range of $\mathbf{X}$; that is,
\begin{equation}
\mathbf{\Pi}_{\ell}^2 = \mathbf{\Pi}_{\ell} = \mathbf{\Pi}_{\ell}^\top.
\end{equation}
Recall that orthogonal projectors on $\mathbb{R}^m$ are non-expansive maps, and that $\mathbf{\Pi}_{\ell}$ satisfies the following optimality property \cite{Quarteroni2016}.
\begin{theorem}\label{thm:proj}
Let $Q = \left\lbrace \mathbf{Q}\in\mathbb{R}^{m\times \ell}~\middle|~\mathbf{Q}^{\top}\mathbf{Q}=\mathbf{I}_{\ell}\right\rbrace$, where $\mathbf{I}_{\ell}\in\mathbb{R}^{\ell\times \ell}$ is the identity matrix, be the set of real matrices whose columns form $\ell$-dimensional orthonormal bases. Then, the projector $\mathbf{\Pi}_{\ell}$ that is defined by \eqref{eq:projector} is the best projector among all projectors of the form $\mathbf{Q}\mathbf{Q}^{\top}$, where $\mathbf{Q}\in Q$, in the sense that
\begin{equation}
\sum_{k=1}^{n}\Vert \mathbf{x}_k - \mathbf{\Pi}_{\ell}\mathbf{x}_k\Vert^{2} = \min_{\mathbf{Q}\in Q} \sum_{k=1}^{n}\Vert\mathbf{x}_k - \mathbf{Q}\mathbf{Q}^{\top}\mathbf{x}_k\Vert^{2} =\sum_{k=r+1}^{\operatorname{rank}(\mathbf{X})}\sigma_{k}^{2}.
\end{equation}
\end{theorem}

Since the singular value decomposition is optimal in the least-squares sense, the proper orthogonal decomposition is suitable for reconstructing the manifold that is formed in $\mathbb{R}^m$ by the columns of the snapshot matrix, when the endomorphism of the underlying system is linear. For nonlinear systems, local applications of the proper orthogonal decomposition on subcollections of the columns of the snapshot matrix may exhibit improved properties \cite{Hediger1990}, since any manifold is locally linear.

\section{Maximum entropy sampling}\label{sec:mess}
Here, we develop the so-called maximum entropy snapshot sampling, which is based on the entropy theory of dynamical systems. Instead of performing the reduction of a basis according to a gauge of information that is based on linear transformations, such as the singular values, we directly reduce the snapshot matrix according to an estimate of the second-order R{\'e}nyi entropy. In particular, our presentation begins with a motivating construction that is based on the invariant second order R{\'e}nyi entropy and on the average sojourn time.

\subsection{Entropy in the discrete setting}
Let $m$ and $n$ be positive integers such that $m\gg n>1$, and consider a system $(\Omega,\mathfrak{A},\mu,f)$. By iterating the endomorphism $f$, we obtain a snapshot matrix
\begin{equation}
\mathbf{X}=(\mathbf{x}_1\mid \mathbf{x}_2\mid\ldots\mid \mathbf{x}_{n})\in\mathbb{R}^{m\times n}.
\end{equation}
Let $\mathbb{I}_n=\{1,2,\ldots,n\}$ be an index set and consider a set $\mathcal{B}=\left\{B_\epsilon(\mathbf{x}_j)\middle|j\in\mathbb{I}_n\right\}$ of open balls 
\begin{equation}
B_\epsilon(\mathbf{x}_j) = \{ \mathbf{y}\in\mathbb{R}^m\mid \Vert \mathbf{x}_j - \mathbf{y}\Vert <\epsilon\}
\end{equation}
that cover $\mathbf{X}$. Then, the second-order R{\'e}nyi entropy of $\mathcal{B}$ is given by
\begin{equation}
H^{2}_{\mu}(\mathcal{B}) = -\log\sum_{j=1}^{n} \mu^2(B_\epsilon(\mathbf{x}_j)) = -\log E_{\mu}(\mu),
\end{equation}
where $E_{\mu}(\mu)$ is the expected value of $\mu$ with respect to $\mu$ itself. When $n$ is sufficiently large, the expected value of $\mu$ can be replaced by the average value of the sequence $(\mu(B_\epsilon(\mathbf{x}_1)),\mu(B_\epsilon(\mathbf{x}_2)),\ldots,\mu(B_\epsilon(\mathbf{x}_n)))$, according to the law of large numbers, while each $\mu(B_\epsilon(\mathbf{x}_j))$ can be approximated by the sample's average sojourn time, according to \eqref{eq:measure}. Under these considerations, we obtain the estimate
\begin{equation}\label{eq:entest1}
\hat{H}^{2}_{\mu}(\mathcal{B}) = -\log\frac{1}{n^2}\sum_{i=1}^{n}\sum_{j=1}^{n} \chi_{B_\epsilon(\mathbf{x}_i)}(\mathbf{x}_j).
\end{equation}
The symmetric matrix $\mathbf{R}_\epsilon\in\{0,1\}^{n\times n}$ whose entries are $\chi_{B_\epsilon(\mathbf{x}_i)}(\mathbf{x}_j)$ is referred to as the recurrence matrix \cite{Marwan2007} that is associated with the sample $\mathbf{X}$ and constitutes a thresholded version of the matrix whose entries are pairwise distances between the columns of $\mathbf{X}$. Entropy estimate \eqref{eq:entest1} can be written as
\begin{equation}\label{eq:entest2}
\eta_\epsilon(\mathbf{X}) = -\log\left(\frac{1}{n^2}\Vert \mathbf{R}_\epsilon\Vert^{2}_\mathrm{F}\right),
\end{equation}
with so-called information potential \cite{Jenssen2010}
\begin{equation}
v_\epsilon(\mathbf{X}) = \frac{1}{n^2}\Vert \mathbf{R}_\epsilon\Vert^{2}_\mathrm{F}\neq 0,
\end{equation}
also known in the literature as the recurrence rate \cite{Marwan2007}. Then, for sufficiently large $n$, the invariant entropy $\mathcal{H}^{2}_\mu(f)$ can be approximated by successive entropy changes in time \cite{Kantz2004}. In the following definition, we collect quantities that are essential for the maximum entropy snapshot sampling.

\begin{definition}
Let $m$ and $n$ be positive integers such that $m\gg n>1$, and consider a snapshot matrix $\mathbf{X}\in\mathbb{R}^{m\times n}$ and -- provided $\epsilon>0$ -- an associated recurrence matrix $\mathbf{R}_\epsilon\in\{0,1\}^{n\times n}$. Then, for each index $j\in\mathbb{I}_n$,
\begin{equation}
\eta_{j}=\eta_\epsilon(\mathbf{X}(:,1:j))=-\log\left(\frac{1}{j^2}\Vert \mathbf{R}_\epsilon(1:j,1:j)\Vert^{2}_\mathrm{F}\right),
\end{equation}
is defined to be the $\epsilon$-Frobenius entropy of $\mathbf{X}(:,1:j)$, with
\begin{equation}
\qquad v_{j}=v_\epsilon(\mathbf{X}(:,1:j)) = \frac{1}{j^2}\Vert \mathbf{R}_\epsilon(1:j,1:j)\Vert^{2}_\mathrm{F}
\end{equation}
being the $\epsilon$-Frobenius potential of $\mathbf{X}(:,1:j)$. Further, 
\begin{equation}
h_{j} = \eta_{j+1} - \eta_{j} = -\log{\frac{v_{j+1}}{v_{j}}}
\end{equation}
is called the $\epsilon$-dynamical entropy of $\mathbf{X}$ at time step $j\in\mathbb{I}_{n-1}$. 
\end{definition}
In the language of dynamical systems and information theory, $h_{j}$ is the information that has been delivered by the iterate $\mathbf{x}_{j+1}=f(\mathbf{x}_j)$, relative to the information that has been delivered by all former iterates.

The $\epsilon$-Frobenius potential is closely related to the so-called correlation sum \cite{Grassberger1983, Kantz2004, Takens1998}, which is the appropriately scaled sum of the entries of the strictly upper triangular part of the recurrence matrix. The correlation sum is used in the definition of the correlation dimension of the manifold that is formed in $\mathbb{R}^m$ by the columns of the snapshot matrix. Entry-wise differential transformations of distance matrices are an essential ingredient of kernel model reduction \cite{Wirtz2014}, and of clustering methods such as the kernel principal component analysis \cite{Scholkopf1998} and the entropy component analysis \cite{Jenssen2010}, in particular, whenever differential optimizations needs to be performed. Further, when the number of available snapshots is sufficiently large, time-delay techniques, which are based on embedding theorems \cite{Takens1981}, enable entropy estimates using a single time-delayed variable. This is not a case of interest in reduced basis methods, since each snapshot requires the solution of a high-dimensional and possibly nonlinear large system of algebraic equations, and hence, we are interested in the smallest possible sample of snapshots that reflects the dynamics of the underlying high-fidelity system.

\subsection{Maximum entropy snapshot sampling}
The key idea of the maximum entropy snapshot sampling is that since the dynamical entropy $h_{j}$ reflects the information gain per iteration step, we require that $h_j>0$ for all $j\in\mathbb{I}_{n-1}$, or equivalently, that the $\epsilon$-Frobenius entropy is a strictly increasing function of the time index $j$. A direct implication of this monotonicity requirement is that the $\epsilon$-Frobenius potential is a strictly decreasing function of the time index.

Let $m$ and $n$ be integers such that $m\gg n>1$. The $\epsilon$-Frobenius potentials $v_{j}$ that are associated with a snapshot matrix $\mathbf{X}\in\mathbb{R}^{m\times n}$ are the instantaneous fractions of all pairs of snapshots $\mathbf{x}_j\in\mathbb{R}^m$ that are closer than a threshold distance $\epsilon$. Since every snapshot $\mathbf{x}_j$ is in its own ball $B_\epsilon(\mathbf{x}_j)$  and at most in the ball of every other snapshot, we have
\begin{equation}
j\leq\Vert\mathbf{R}(1:j,1:j)\Vert_\mathrm{F}^2\leq j^2\Leftrightarrow 1/j \leq v_j\leq 1\qquad\forall j\in\mathbb{I}_{n},
\end{equation}
where we dropped the $\epsilon$ subscript from the recurrence matrix for notational simplicity. The following theorem transfers the monotonicity of the $\epsilon$-Frobenius entropy to the structure of the recurrence matrix. In particular, \cref{thm:ident} states that if the $\epsilon$-Frobenius entropy attains its maximum value $\log{j}$ at each timestep $j$ (or if the $\epsilon$-Frobenius potential attains its minimum value $1/j$), then $\mathbf{R}$ is the identity matrix, and hence, each snapshot is only in its own ball.

\begin{theorem}\label{thm:ident}
Let $m$ and $n$ be integers such that $m\gg n>1$, and consider a matrix $\mathbf{X}$ with columns $\mathbf{x}_j\in\mathbb{R}^m$ for all $j\in\mathbb{I}_n$. Provided $\epsilon>0$, the $\epsilon$-Frobenius entropy \eqref{eq:entest2} is a strictly increasing function of the time index, if, and only if, $\Vert \mathbf{x}_j - \mathbf{x}_k\Vert \geq \epsilon$ for all distinct $j$, $k\in\mathbb{I}_n$.
\begin{proof}
Assume that the $\epsilon$-Frobenius entropy is a strictly increasing function of the time index; that is, $\eta_{j+1}>\eta_{j}$ for all $j\in\mathbb{I}_{n-1}$. Then, due to the monotonicity of the $\log$ function, we obtain $v_{j+1}<v_{j}$ for all $j\in\mathbb{I}_{n-1}$. By its definition, the information potential satisfies the initial value problem 
\begin{equation}\label{eq:recprob}
v_1 = 1,\qquad v_{j+1} = \frac{j^2}{(j+1)^2} v_j + \frac{1}{(j+1)^2}\delta_{j}\quad\forall j\in\mathbb{I}_{n-1},
\end{equation}
where, due to the symmetry of each block $\mathbf{R}(1:j,1:j)$,
\begin{equation}\label{eq:delta}
\delta_{j} = 2\sum_{k=1}^{j}\mathbf{R}(j+1,k) + 1
\end{equation}
is an odd number that is bounded as $1\leq\delta_{j}\leq 2j+1$, see \cref{fig:proof}.
\begin{figure}
  \centering
  \includegraphics[width=0.3\textwidth]{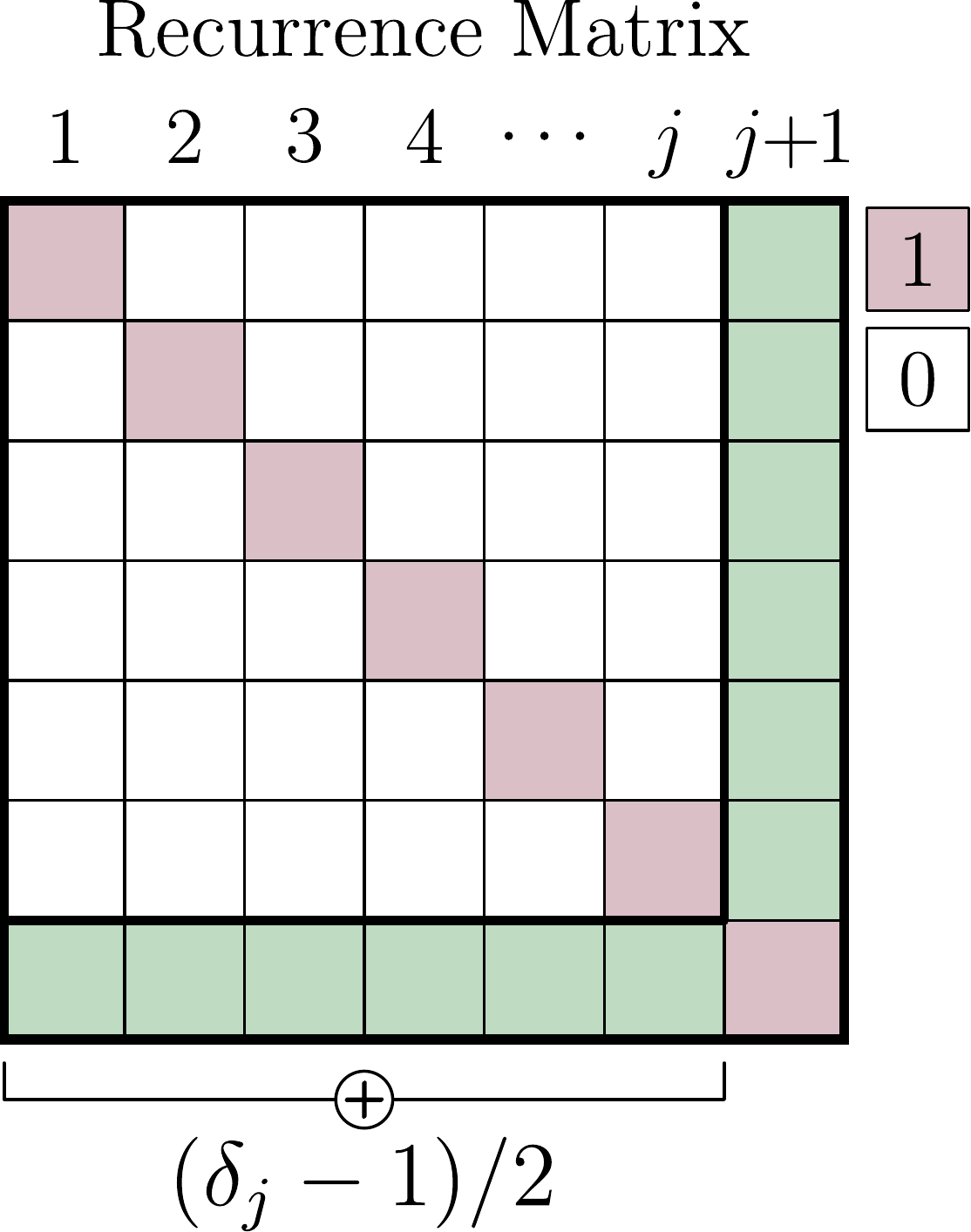}\label{fig:proof}
  \caption{The construction that is used in the proof of \cref{thm:ident}.}
\end{figure}
Thus,
\begin{equation}\label{eq:diffpot}
v_{j+1}-v_{j} = \frac{\delta_{j}-(2j+1)v_{j}}{(j+1)^2}
\end{equation}
and since $v_{j+1}<v_{j}$, we obtain
\begin{equation}\label{eq:ineq}
\delta_{j}<(2j+1)v_{j}.
\end{equation}
By recursively evaluating \eqref{eq:ineq} for $j\in\mathbb{I}_{n-1}$, we obtain $v_{j}=1/j$ and $\delta_{j} < 3$ for all $j\in\mathbb{I}_{n-1}$. Since $\delta_{j}\geq 1$ is an odd number, we conclude that $\delta_{j}=1$ for all $j\in\mathbb{I}_{n-1}$. Equation \eqref{eq:delta} entails that $\mathbf{R}(j+1,k)=0$ for all $k\in\mathbb{I}_j$, and hence, $\mathbf{R}=\mathbf{I}_n$. By the definition of the entries of $\mathbf{R}$, we obtain $\Vert \mathbf{x}_j - \mathbf{x}_k\Vert \geq \epsilon$ for all distinct $j$, $k$. The converse implication is evident and thus, the proof is complete.
\end{proof}
\end{theorem}

\cref{thm:ident} suggests that to generate a reduced basis for the range of $\mathbf{X}$, we only need the columns of $\mathbf{X}$ whose $\epsilon$-Frobenius entropy is an increasing function of the column index, while it provides a memory efficient way for sampling the required columns through initial value problem \eqref{eq:recprob}. The following theorem provides an error bound between the original and the back-projected columns of $\mathbf{X}$, when a basis is generated from non-recurrent columns of $\mathbf{X}$. 
\begin{theorem}\label{thm:main}
Let $\ell$, $m$, $n$ be integers such that $m\gg n>1$ and $\ell\leq n$, and consider a matrix $\mathbf{X}\in\mathbb{R}^{m\times n}$ whose columns are $\mathbf{x}_j\in\mathbb{R}^m$, where $j\in\mathbb{I}_n$. Provided $\epsilon>0$, form the matrix $\mathbf{Y}_\ell\in\mathbb{R}^{m\times\ell}$ whose columns $\mathbf{y}_{k}$, with $k\in\mathbb{I}_{\ell}$, constitute a subset of the columns of $\mathbf{X}$ for which the $\epsilon$-Frobenius entropy estimate \eqref{eq:entest2} is a strictly increasing function of the column index $j$. Then, for any matrix $\mathbf{Q}\in\mathbb{R}^{m\times\ell}$ whose columns form an orthonormal basis for the range of $\mathbf{Y}_{\ell}$, we have that each $\hat{\mathbf{x}}_j=\mathbf{Q}\mathbf{Q}^{\top}\mathbf{x}_j$ satisfies the error bound $\Vert \mathbf{x}_j - \hat{\mathbf{x}}_j\Vert < \epsilon$ for all $j\in\mathbb{I}_n$.
\begin{proof}
When a high-fidelity snapshot $\mathbf{x}_{j}$ is a column of $\mathbf{Y}_{\ell}$, then there exists some $\mathbf{r}\in\mathbb{R}^{\ell}$, such that $\mathbf{x}_j=\mathbf{Q}\mathbf{r}$, and hence,
\begin{equation}
\hat{\mathbf{x}}_j = \mathbf{Q}\mathbf{Q}^{\top}\mathbf{x}_j =  \mathbf{Q}\mathbf{Q}^{\top}\mathbf{Q}\mathbf{r} = \mathbf{Q}\mathbf{r} = \mathbf{x}_j,
\end{equation}
since $\mathbf{Q}^{\top}\mathbf{Q}=\mathbf{I}_\ell$. If $\mathbf{x}_{j}$ is not a column of $\mathbf{Y}_{\ell}$, then, since $\mathbf{y}_{k}=\mathbf{Q}\mathbf{Q}^{\top}\mathbf{y}_{k}$ for any column $\mathbf{y}_{k}$ of $\mathbf{Y}_{\ell}$, we have
\begin{align}
\Vert\mathbf{x}_j - \hat{\mathbf{x}}_j\Vert & = \Vert \mathbf{x}_j - \mathbf{y}_k + \mathbf{y}_k - \hat{\mathbf{x}}_j\Vert\\
& = \Vert \mathbf{x}_j - \mathbf{y}_k + \mathbf{Q}\mathbf{Q}^\top\mathbf{y}_k - \mathbf{Q}\mathbf{Q}^\top\mathbf{x}_j\Vert\\
& = \Vert (\mathbf{I}_m - \mathbf{Q}\mathbf{Q}^\top)(\mathbf{x}_j-\mathbf{y}_k)\Vert\\
& \leq \Vert \mathbf{x}_j-\mathbf{y}_k\Vert < \epsilon,
\end{align}
where the first inequality follows from the fact that $\mathbf{I}_m-\mathbf{Q}\mathbf{Q}^{\top}$ is an orthogonal projector, and as such, it is a non-expansive map, while the last inequality follows from the maximum entropy snapshot sampling, since for each $j\in\mathbb{I}_{n}$ there exists some $k\in\mathbb{I}_{\ell}$ such that $\mathbf{x}_j\in B_{\epsilon}(\mathbf{y}_{k})$.
\end{proof}
\end{theorem}

If the columns of $\mathbf{X}$ are the snapshots of a stationary dynamical system, then the limit $v_{\ast}$ of the Frobenius potential exists, according to the presentation in \cref{sec:dynsys}. When this is the case, the following theorem suggests that the sampling process can be interrupted at a particular point of time, and guarantees that the basis generated in \cref{thm:main} keeps the error between the exact and the back-projected future snapshots bounded up to a given time-horizon.
\begin{theorem}\label{thm:lim}
If there exists an integer $j_{\ast}$ such that $v_{j+1}=v_{j}=v_{\ast}$ for all $j\geq j_{\ast}$ and an index $i\geq 1$ that satisfies the horizon inequality
\begin{equation}\label{eq:horizon}
i(1-v_{\ast}) < j_{\ast}v_{\ast} - \frac{1}{2}(1-v_{\ast}),
\end{equation}
then each $\hat{\mathbf{x}}_{j_{\ast}+i}=\mathbf{Q}\mathbf{Q}^{\top}\mathbf{x}_{j_{\ast}+i}$, where $\mathbf{Q}$ is the reduced basis in \cref{thm:main}, satisfies the error estimate $\Vert \mathbf{x}_{j_{\ast}+i}-\hat{\mathbf{x}}_{j_{\ast}+i}\Vert<4\epsilon$.
\begin{proof}
Since $v_{j+1}=v_{j}=v_{\ast}$ for all $j\geq j_{\ast}$, from \eqref{eq:recprob} we obtain $\delta_{j}=(2j+1)v_{\ast}$ or $\delta_{j_{\ast}+i}=(2j_{\ast}+2i+1)v_{\ast}$ for all integers $i\geq 0$. Further, since
\begin{align}
i(1-v_{\ast}) < j_{\ast}v_{\ast} - \frac{1}{2}(1-v_{\ast}) & \Leftrightarrow (2j_{\ast}+2i+1)v_{\ast}>2i+1\\
                                                           & \Leftrightarrow \delta_{j_{\ast}+i} > 2i+1,
\end{align}
each snapshot $\mathbf{x}_{j_{\ast}+i}$ is in the ball of some snapshot $\mathbf{x}_{j}$ with $j\in\mathbb{I}_{j_{\ast}}$. If $\mathbf{x}_{j_\ast+i}$ is in the ball of a snapshot $\mathbf{y}_k$ of the sampled matrix $\mathbf{Y}_{\ell}$ in \cref{thm:main}, then $\Vert\mathbf{x}_{j_{\ast}+i}-\mathbf{y}_k\Vert<\epsilon$, while otherwise, an application of the triangle inequality results in $\Vert\mathbf{x}_{j_{\ast}+i}-\mathbf{y}_k\Vert < 2\epsilon$. Thus, 
\begin{align}
\Vert \mathbf{x}_{j_{\ast}+i} - \hat{\mathbf{x}}_{j_{\ast}+i}\Vert &= \Vert \mathbf{x}_{j_{\ast}+i} -\mathbf{y}_k +\mathbf{y}_k - \hat{\mathbf{x}}_{j_{\ast}+i}\Vert\\
                           & \leq  \Vert \mathbf{x}_{j_{\ast}+i} -\mathbf{y}_k\Vert +\Vert \mathbf{y}_k - \hat{\mathbf{x}}_{j_{\ast}+i}\Vert\\
                           & < 2\epsilon + \Vert\mathbf{Q}\mathbf{Q}^{\top}\mathbf{y}_{k} - \mathbf{Q}\mathbf{Q}^{\top}\mathbf{x}_{j_{\ast}+i}\Vert\\
                           & \leq 2\epsilon + \Vert\mathbf{y}_{k} - \mathbf{x}_{j_{\ast}+i}\Vert\\
                           & < 4\epsilon,
\end{align}
where the third inequality follows from the fact that $\mathbf{Q}\mathbf{Q}^{\top}$ is an orthogonal projector, and as such, it is a non-expansive map.
\end{proof}
\end{theorem}

\cref{thm:lim} suggests a criterion for interrupting the sampling process from the high-fidelity solver. In particular, the sampling process can be interrupted when $|v_{j+1}-v_j|<\varepsilon$ for some tolerance $0<\varepsilon\ll 1$, so that the order $\epsilon$ error bound that is satisfied by the maximum entropy snapshot sampling basis is approximately preserved for the future data, up to an index $i$ that satisfies the horizon inequality \eqref{eq:horizon}. A similar criterion for the $\epsilon$-dynamical entropy may be also used, that is, $|h_{j+1}-h_{j}|<\varepsilon$.


In the remaining part of this work, we present numerical experiments that validate the maximum entropy snapshot sampling as an efficient compression and reduction method.

\section{Numerical experiments}
We test the error bound in \cref{thm:main} and the performance of the maximum entropy snapshot sampling (MESS) by obtaining a reduced basis with the economy-size orthogonal-triangular decomposition QR of the MESS reduced samples. The reference method is the singular value decomposition (SVD), as implemented in the $\text{MATLAB}^\text{\textregistered}$ routines \texttt{svd}, with economy option, or \texttt{svds}, with the number of basis vectors being the same with the MESS reduced samples.

\begin{figure}
  \centering
  \includegraphics[width=0.9\textwidth]{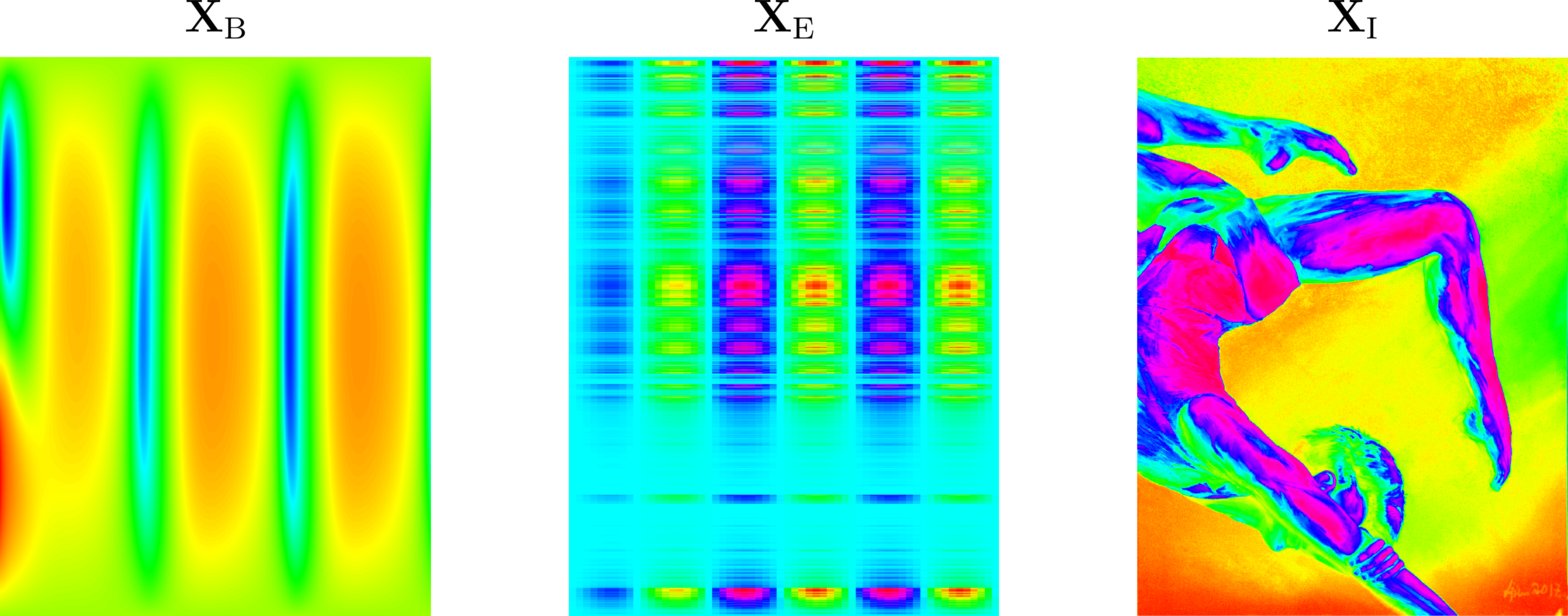}
  \caption{The data-sets $\mathbf{X}_\mathrm{B}$, $\mathbf{X}_\mathrm{E}$, and $\mathbf{X}_\mathrm{I}$ that are used for verifying the error bount in \cref{thm:main} and the performance of the MESS. The data are normalized in the intervals $[0,1]$, $[-1,1]$, $[0,1]$, respectively.}\label{fig:snaps}
\end{figure}

For the first part of our verification tests we use three sets of data that are arranged in the matrices $\mathbf{X}_\mathrm{B}$, $\mathbf{X}_\mathrm{E}$, $\mathbf{X}_\mathrm{I}\in\mathbb{R}^{m\times n}$, with
\begin{equation}
(m,n)\in\{(200\,000,2\,001),(1\,043\,571,60),(1\,687,1\,181)\}, 
\end{equation}
respectively, see \cref{tab:data}. The data-set $\mathbf{X}_\mathrm{B}$ is a snapshot matrix obtained from the so-called brusselator system \cite{Tyson1973}, which models diffusion in chemical reaction, using a modified version of the $\text{MATLAB}^\text{\textregistered}$ solver \texttt{brussode}. At each one of the $n=2\,001$ snapshots, the state of the brusselator system is specified by determining the concentrations of $m=200\,000$ reacting species.
\begin{table}[tbhp]
{\footnotesize\caption{Information about the source problems of the data that are used for the verification tests.}\label{tab:data}
\begin{center}
\begin{tabular}{|c|c|c|c|c|c|c|}\hline
\bf Model & \bf Data &\bf Class & \bf Solution method\\ [0.5ex] 
\hline
Brusselator & $\mathbf{X}_\mathrm{B}$ & Deterministic, nonlinear IVP & \texttt{ode15s}\\
Electro-quasistatic & $\mathbf{X}_\mathrm{E}$ & Deterministic, nonlinear IBVP & FEM/BDF(1)\\
Digital image & $\mathbf{X}_\mathrm{I}$ & Non-deterministic & N/A\\
\hline
\end{tabular}
\end{center}
}
\end{table}
The data-set $\mathbf{X}_\mathrm{E}$ is a snapshot matrix obtained with an in-house finite element with backward differentiation in time solver (FEM/BDF(1)) of a transient electro-quasistatic initial-boundary value problem \cite{Haus1989}, implemented in FREEFEM. The computational domain partly consists of strongly nonlinear field grading material \cite{Christen2010}, while a time-harmonic excitation is gradually introduced, within one period $0.02\,\mathrm{s}$. In particular, at each one of the $n=60$ snapshots (three periods), the state of the electro-quasistatic system is specified by determining the scalar electric potential at $m=1\,043\,571$ nodes in the computational domain. The third data-set $\mathbf{X}_\mathrm{I}$ is a digital image of a painting, used with the written permission of the painter \cite{Lipken}. Normalized versions of these data are depicted in \cref{fig:snaps}. Although the data-set $\mathbf{X}_\mathrm{I}$ is not deterministic and the notion of time is not defined, there is no argument that prohibits us from using the maximum entropy snapshot sampling as a compression algorithm, since determinism is not a requirement for the validity of \cref{thm:main}. 

\begin{figure}
  \centering
  \includegraphics[width=0.99\textwidth]{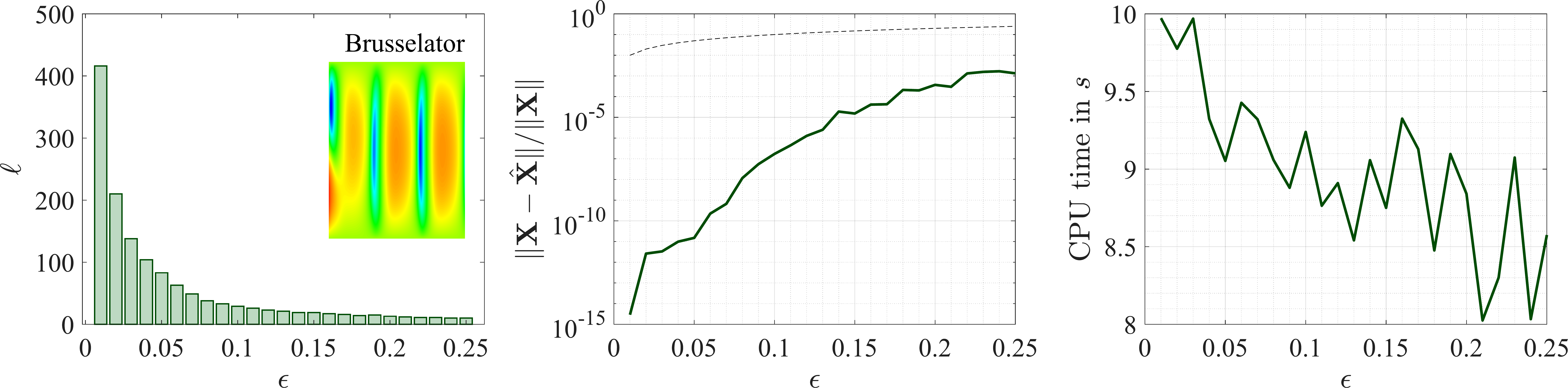}
  \caption{The number of reduced basis vectors (left), the relative Euclidean error (middle), and the CPU time that is required to obtain a reduced basis (right) with MESS versus $\epsilon$ for $\mathbf{X}=\mathbf{X}_\mathrm{B}$.}\label{fig:XBres}
\end{figure}

\begin{figure}
  \centering
  \includegraphics[width=0.99\textwidth]{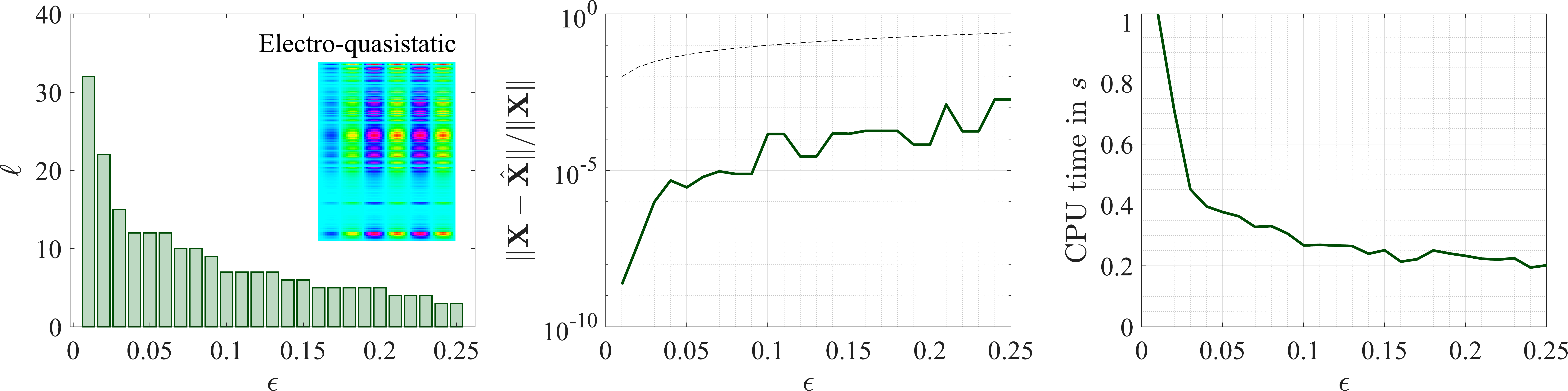}
  \caption{The number of reduced basis vectors (left), the relative Euclidean error (middle), and the CPU time that is required to obtain a reduced basis (right) with MESS versus $\epsilon$ for $\mathbf{X}=\mathbf{X}_\mathrm{E}$.}\label{fig:XEres}
\end{figure}

\begin{figure}
  \centering
  \includegraphics[width=0.99\textwidth]{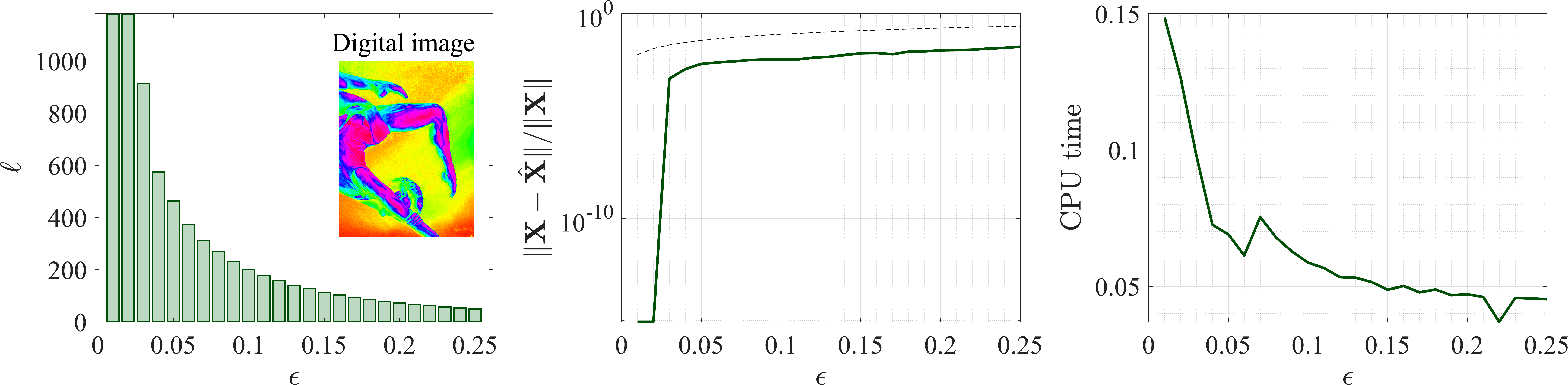}
  \caption{The number of reduced basis vectors (left), the relative Euclidean error (middle), and the CPU time that is required to obtain a reduced basis (right) with MESS versus $\epsilon$ for $\mathbf{X}=\mathbf{X}_\mathrm{I}$.}\label{fig:XIres}
\end{figure}

Each matrix $\mathbf{X}\in\{\mathbf{X}_\mathrm{B},\mathbf{X}_\mathrm{E},\mathbf{X}_\mathrm{I}\}$, is compressed with the MESS framework for $\epsilon\in\{0.01,0.02,\ldots,0.25\}$, as a percentage of the maximum pairwise distance of its columns. In \cref{fig:XBres}, \cref{fig:XEres}, and \cref{fig:XIres} we plot the number of reduced basis vectors (left), the relative Euclidean error (middle), and the CPU time (in seconds) that is required to obtain a reduced basis (right) versus $\epsilon$ for each data-set, respectively. On the left of these figures, we observe that the user parameter $\epsilon$ determines the degree of reduction, with smaller values of $\epsilon$ resulting in less reduction, since smaller balls yield a smaller number of recurrent states. In the error plots of the same figures, the dashed lines correspond to the function with values $f(\epsilon)=\epsilon$, supporting the validity of the derived error estimate, independent of the nature of the source of these data. It is worth mentioning that for less regular data, that is for $\mathbf{X}_\mathrm{I}$, the MESS error is closer to the upper bound $\epsilon$, suggesting that the additional regularity of the deterministic data improves the MESS reduction framework. On the right plots of the same figures, the CPU time is the sum of the sampling time and the time needed to perform an economy-size orthogonal-triangular decomposition QR, without permutations, using the standard \texttt{qr} routine from $\text{MATLAB}^\text{\textregistered}$.

\begin{figure}
  \centering
  \includegraphics[width=0.85\textwidth]{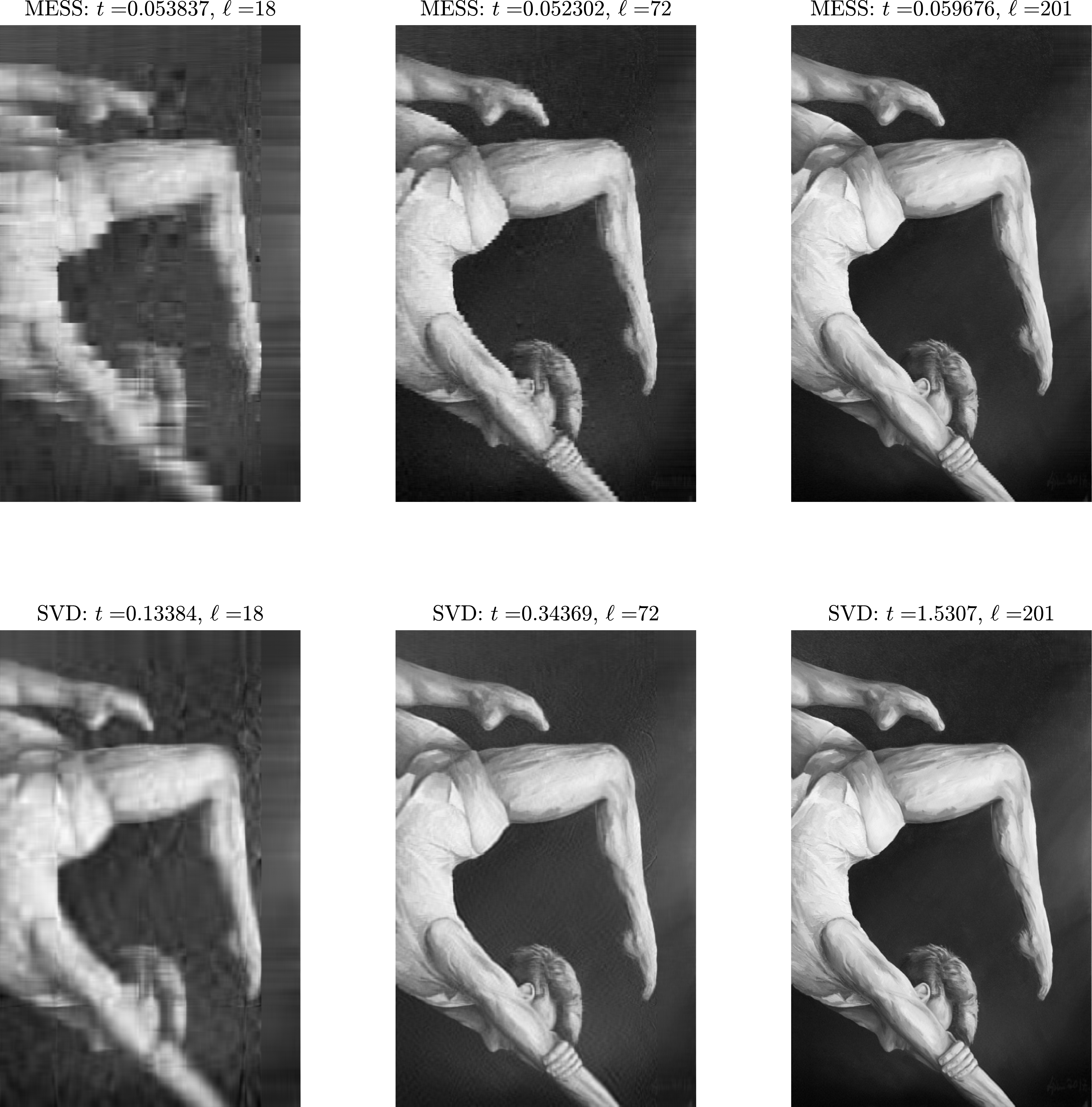}
  \caption{MESS (upper row) and SVD (lower row) based compression of the digital image $\mathbf{X}_\mathrm{I}$, with $m=1687$ and $n=1181$, for $\epsilon\in\{0.4,0.2,0.1\}$ (from left to right). The size of each SVD generated basis is chosen to match the size of the MESS generated length. In terms of CPU time, the MESS outperforms the SVD in all cases.}\label{fig:imag}
\end{figure}

To visualize the result of the MESS basis and to examine its performance against the SVD, we work with the matrix $\mathbf{X}_\mathrm{I}$. In \cref{fig:imag}, we report the CPU times for MESS and SVD generated bases of the same length, and we observe that the MESS outperforms SVD, in the examined cases. When attempting to match the error between the MESS and SVD compression, $\Vert\mathbf{X}-\hat{\mathbf{X}}_{\mathrm{MESS}}\Vert/\Vert\mathbf{X}-\hat{\mathbf{X}}_{\mathrm{SVD}}\Vert\cong 1$, then for $\ell_\mathrm{SVD}\in\{25,50,100\}$ SVD basis vectors we obtain $\ell_\mathrm{MESS}\in\{49,85,292\}$ MESS basis vectors, entailing that the MESS provides a smaller degree of reduction than SVD, but still at a fraction of the time, with $t_\mathrm{MESS}/t_\mathrm{SVD}\in\{0.42,0.25,0.16\}$, respectively. The timing results in \cref{tab:timing}, confirm that MESS outperforms the SVD also for the data $\mathbf{X}_\mathrm{B}$, $\mathbf{X}_\mathrm{E}$ for $\epsilon\in\{0.01,0.05\}$, while the effect of calling the basis generation factorization routine before and after the MESS, in particular, for mid- and large-scale problems with many snapshots, becomes apparent.
\begin{remark}
For the matching accuracy experiments, we chose arbitrary numbers of SVD basis vectors. In terms of performance, the situation for SVD is actually worse than the one presented here, since the computational cost of criterion \eqref{eq:svdinfo} has to be taken into account. 
\end{remark}
\begin{table}[tbhp]
{\footnotesize\caption{Timing information for the brusselator and the electric quasistatic snapshot matrices. The \texttt{svd} routine with the economy option has been used for these times, since it performed better than \texttt{svds}.}\label{tab:timing}
\begin{center}
\begin{tabular}{|c|c|c|c|}\hline
\bf Model & \bf MESS for $\boldsymbol{\epsilon=0.01}$ & \bf MESS for $\boldsymbol{\epsilon=0.05}$ & \bf SVD\\
\hline
Brusselator & $10.7~\mathrm{s}$ & $8.4~\mathrm{s}$ & $53~\mathrm{s}$\\
Electro-quasistatic & $1.1~\mathrm{s}$ & $0.36~\mathrm{s}$ & $1.1~\mathrm{s}$\\
\hline
\end{tabular}
\end{center}
}
\end{table}

\begin{figure}
  \centering
  \includegraphics[width=0.5\textwidth]{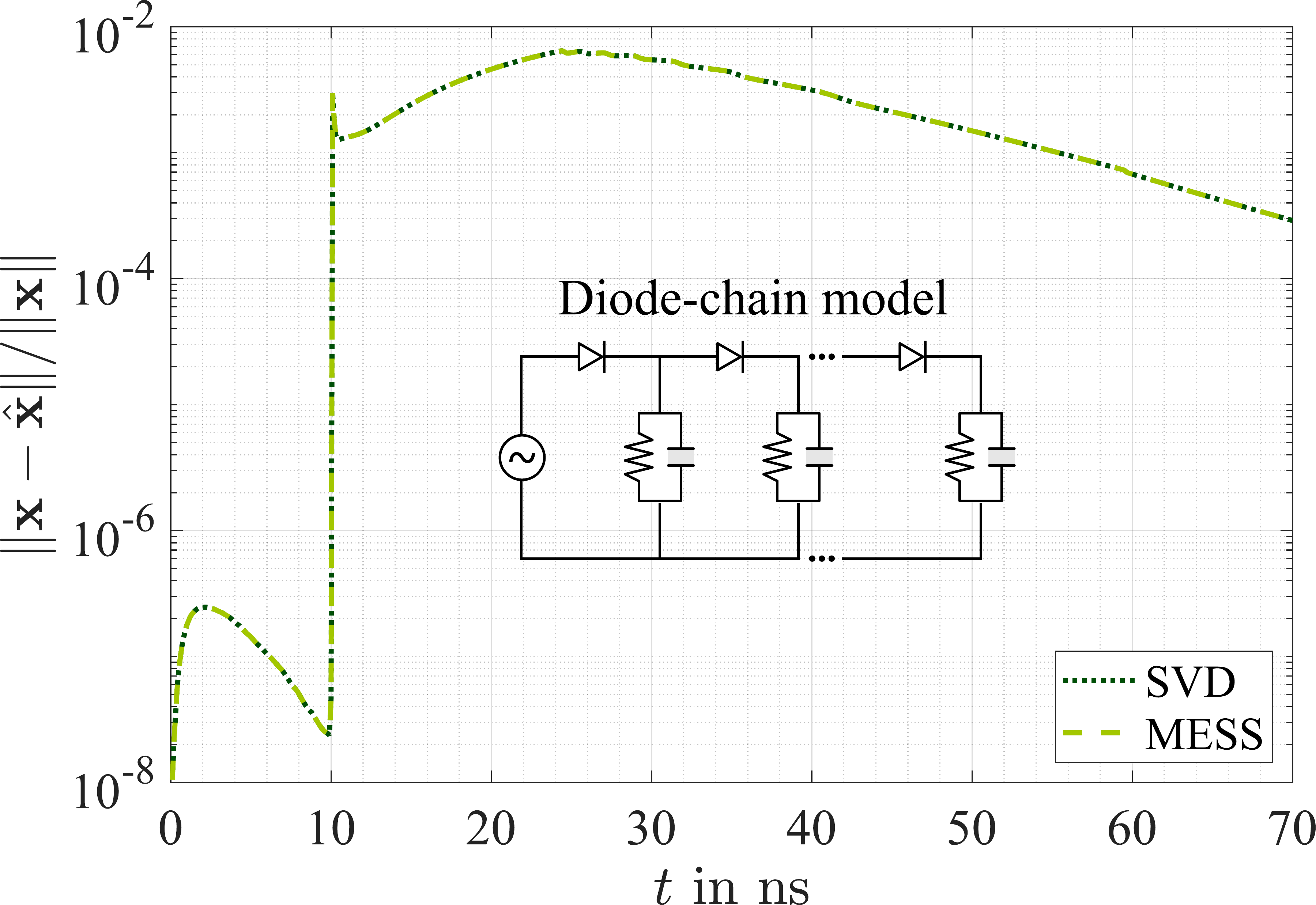}
  \caption{The relative Euclidean error of the MESS and SVD reduced basis solutions with respect to the high-fidelity solution of a diode chain model. Here, $\mathbf{x}$ is the instantaneous reference solution, while $\hat{\mathbf{x}}$ represents the reduced basis solutions for $\epsilon=0.0325$ and $\ell=31$ basis vectors for both MESS and SVD.}\label{fig:rbmor}
\end{figure}

As a last experiment, the MESS framework is used for the reduction of a nonlinear differential-algebraic diode-chain model \cite{Verhoeven2007}, with $m = 40002$ degrees of freedom. We employ the BDF(1) scheme to obtain $n = 701$ high-fidelity snapshots in the interval $[0,70]~\mathrm{ns}$. The MESS and SVD bases are obtained from the resulting snapshot matrix $\mathbf{X}\in\mathbb{R}^{m \times n}$. The number of SVD basis vectors is taken to be equal to the number of MESS basis vectors, while in the Newton iterations, least-squares approximations of the Jacobian matrix are used. For $\epsilon=0.0325$ we obtain $\ell_\mathrm{MESS}=\ell_\mathrm{SVD} = 31$ basis vectors, while the offline stage of the SVD requires $8.9$ times the CPU time that is required by the MESS, for identical accuracy, as depicted in \cref{fig:rbmor}. It is worth mentioning that the $\epsilon$ value has been chosen so that both MESS and SVD result in numerically stable reduced models. 

\section{Conclusions}
We introduced a general and efficient sampling method that is appropriate for compression and reduced basis model reduction. We proved that the maximum entropy snapshot sampling (MESS) satisfies a problem-independent error bound, while it enables the usage of arbitrary basis generation algorithms. In the case of deterministic evolution problems, we proved that the reliability of the maximum entropy snapshot sampling bases, in terms of error control, depends on the recurrence properties of the source system and we delivered the so-called horizon inequality, which guarantees the derived error bound. In terms of computational complexity, we numerically verified that our implementation outperforms the singular value decomposition, although a higher degree of reduction is enabled by the latter. An important benefit of MESS over the singular value decomposition is that it is memory-friendly and can be used in parallel with the high-fidelity solver. This last feature enables the reduction of large-scale problems for which the basis generation stage becomes a bottleneck, when using the proper orthogonal decomposition, if feasible at all. Hence, MESS can be used for very large-scale problems that result from finite element discretizations of three-dimensional problems, without overloading the offline stage. Further, since the MESS relies on pairwise distance computations, it is possible to employ CPU/GPU parallelization to further speed up the offline stage.

\section*{Acknowledgements}
The authors acknowledge the contribution of Marcus Bannenberg and Michael Günther (IMACM, Chair of Applied Mathematics, University of Wuppertal) in the reduced basis model reduction results for the diode chain model.

\bibliographystyle{siamplain}
\bibliography{references}
\end{document}